\documentclass{amsart}
\usepackage{amssymb,amsmath, amsthm,latexsym}
\usepackage{graphics}
\usepackage{amscd}
\usepackage{graphics}
\usepackage{here}
\newcommand{\cal}[1]{\mathcal{#1}}
\theoremstyle{plain}

\parskip=\bigskipamount

\let\egthree=\phi
\let\phi=\varphi
\let\varphi=\egthree

\begin{document}
\title{Train tracks and the Gromov boundary of the
complex of curves}
\author{Ursula Hamenst\"adt}
\thanks
{Partially supported by Sonderforschungsbereich 611}
\date{January 28, 2005}

\maketitle

\section{Introduction}

Consider a compact oriented surface $S$ of
genus $g\geq 0$ from which $m\geq 0$ points,
so-called \emph{punctures}, have been deleted.
We require that $3g-3+m\geq 2$; this rules out a
sphere with at most 4 punctures and a torus with
at most one puncture.

In [Ha], Harvey defined the \emph{complex of curves} ${\cal C}(S)$
for $S$. The vertices of this complex are free homotopy classes of
simple closed curves on $S$. The simplices in ${\cal C}(S)$ are
spanned by collections of such curves which can be realized
disjointly. Thus the dimension of ${\cal C}(S)$ equals
$3g-3+m-1$ (recall that $3g-3+m$ is the number of curves in a
\emph{pants decomposition} of $S$).

The \emph{extended mapping class group} $\tilde {\cal M}_{g,m}$ of
$S$ is the group of all isotopy classes of homeomorphisms of $S$.
It acts naturally on the complex of curves as a group of
simplicial automorphisms. Even more is true: If $S$ is not a torus
with 2 punctures then the extended mapping class group is
\emph{precisely} the group of simplicial automorphisms of ${\cal
C}(S)$ (see [I] for references and a sketch of the proof).

Providing each simplex in ${\cal C}(S)$ with the standard
euclidean metric of side-length 1 equips the complex of curves
with the structure of a geodesic metric space whose isometry group
is just $\tilde {\cal M}_{g,m}$ (except for the twice punctured
torus). However, this metric space is not locally compact. Masur
and Minsky [MM1] showed that nevertheless the geometry of ${\cal
C}(S)$ can be understood quite explicitly. Namely, ${\cal C}(S)$
is hyperbolic of infinite diameter. Recall that for some $\delta
>0$ a geodesic metric space is \emph{$\delta$-hyperbolic in the
sense of Gromov} if it satisfies the \emph{$\delta$-thin triangle
condition}: For every geodesic triangle with sides $a,b,c$ the
side $c$ is contained in the $\delta$-neighborhood of $a\cup b$.
Later Bowditch [B] gave a simplified proof of the result of Masur
and Minsky which can also be used to compute explicit bounds for
the hyperbolicity constant $\delta$.

A $\delta$-hyperbolic geodesic metric space $X$ admits a
\emph{Gromov boundary} which is defined as follows. Fix a point
$p\in X$ and for two points $x,y\in X$ define the \emph{Gromov
product} $(x,y)_p=\frac{1}{2}(d(x,p)+d(y,p)-d(x,y))$. Call a
sequence $(x_i)\subset X$ \emph{admissible} if $(x_i,x_j)_p\to
\infty$ $(i,j\to \infty)$. We define two admissible sequences
$(x_i),(y_i)\subset X$ to be \emph{equivalent} if $(x_i,y_i)_p\to
\infty$. Since $X$ is hyperbolic, this defines indeed an
equivalence relation (see [BH]). The Gromov boundary $\partial X$
of $X$ is then the set of equivalence classes of admissible
sequences $(x_i)\subset X$. It carries a natural Hausdorff
topology with the property that the isometry group of $X$ acts on
$\partial X$ as a group of homeomorphisms. For the complex of
curves, the Gromov boundary was determined by Klarreich [K].

For the formulation of Klarreich's result, recall that a
\emph{geodesic lamination} for a complete hyperbolic structure of
finite volume on $S$ is a \emph{compact} subset of $S$ which is
foliated into simple geodesics. A simple closed geodesic on $S$ is
a geodesic lamination with a single leaf. The space ${\cal L}$ of
geodesic laminations on $S$ can be equipped with the Hausdorff
topology for compact subsets of $S$. With respect to this
topology, ${\cal L}$ is compact and metrizable. A geodesic
lamination is called \emph{minimal} if each of its half-leaves is
dense. A minimal geodesic lamination $\lambda$ \emph{fills up $S$}
if every simple closed geodesic on $S$ intersects $\lambda$
transversely, i.e. if every complementary component of $\lambda$
is an ideal polygon or a once punctured ideal polygon with
geodesic boundary [CEG].

A geodesic lamination is \emph{maximal} if its complementary
regions are all ideal triangles or once punctured monogons. Note
that a geodesic lamination can be both minimal and maximal (this
unfortunate terminology is by now standard in the literature).
Each geodesic lamination $\lambda$ is a \emph{sublamination} of a
maximal lamination, i.e. there is a maximal lamination which
contains $\lambda$ as a closed subset [CEG]. For any minimal geodesic
lamination $\lambda$ which fills up $S$, the number of geodesic
laminations $\mu$ which contain $\lambda$ as a sublamination is
bounded by a universal constant only depending on the topological
type of the surface $S$. Namely, each such lamination $\mu$ can be
obtained from $\lambda$ by successively subdividing complementary
components $P$ of $\lambda$ which are different from an ideal triangle
or a once punctured monogon by adding a simple geodesic line which either
connects two non-adjacent cusps of $P$ or
goes around a puncture in the interior of $P$. Notice that every
leaf of $\mu$ which is not contained in $\lambda$ is necessarily isolated
in $\mu$. 

We say that a sequence $(\lambda_i)\subset {\cal L}$
\emph{converges in the coarse Hausdorff topology} to a
minimal lamination $\mu$ which fills up $S$ if
every accumulation point of $(\lambda_i)$ with respect
to the Hausdorff topology contains $\mu$ as a sublamination.
We equip the space ${\cal B}$ of minimal geodesic
laminations which fill up $S$ with the following topology. A set
$A\subset {\cal B}$ is closed if and only if
for every sequence
$(\lambda_i)\subset A$ which converges in the coarse
Hausdorff topology to a lamination $\lambda\in {\cal B}$
we have $\lambda\in A$.
We call this topology on
${\cal B}$ the \emph{coarse Hausdorff topology}.
Using this terminology, Klarreich's
result [K] can be formulated as follows.

\bigskip

{\bf Theorem:} {\it \begin{enumerate}
\item There is a natural homeomorphism $\Lambda$ of ${\cal B}$
equipped with the coarse Hausdorff topology onto the Gromov
boundary $\partial {\cal C}(S)$ of the complex of curves
${\cal C}(S)$ for $S$.
\item For $\mu\in {\cal B}$
a sequence $(c_i)\subset {\cal C}(S)$ is admissible and
defines the point $\Lambda(\mu)\in \partial {\cal C}(S)$ if and
only if $(c_i)$ converges in the coarse Hausdorff
topology to $\mu$.
\end{enumerate}
}

\bigskip

In the paper [K], Klarreich formulates her result using
\emph{measured foliations} on the surface $S$, i.e. topological
foliations $F$ on $S$ equipped with a \emph{transverse translation
invariant measure}. The space ${\cal M\cal F}$ of measured
foliations can be equipped with the weak$^*$-topology which is
metrizable and hence Hausdorff. This topology projects to a
metrizable topology on the space ${\cal P\cal M\cal F}$ of
projective measured foliations which is the quotient of ${\cal
M\cal F}$ under the natural action of the positive half-line
$(0,\infty)$. A topological foliation on $S$ is called
\emph{minimal} if it does not contain a trajectory which is a
simple closed curve. For every minimal topological foliation $F$,
the set of projective measured foliations whose support equals $F$
is a closed subset of ${\cal P\cal M\cal F}$. It follows that the
quotient ${\cal Q}$ of the space of
minimal projective measured foliations under the measure
forgetting equivalence relation is a Hausdorff space as well.
Note that the extended mapping class group of $S$ acts on ${\cal
Q}$ as a group of homeomorphisms. Klarreich shows that ${\cal Q}$
can be identified with the Gromov boundary of the
complex of curves.

There is a natural map $\iota$ which assigns to a measured
foliation $F$ on $S$ a \emph{measured geodesic lamination}
$\iota(F)$, i.e. a geodesic lamination $\lambda$ together with a
transverse translation invariant measure supported in $\lambda$.
The geodesic lamination $\lambda$ is the closure of the set of
geodesics which are obtained by straightening the non-singular
trajectories of the foliation (see [L] for details), together with
the natural image of the transverse measure. A measured geodesic
lamination can be viewed as a locally finite Borel measure on the
space of unoriented geodesics in the hyperbolic plane which is
invariant under the action of the fundamental group of $S$. Thus
the space ${\cal M\cal L}$ of measured geodesic laminations on $S$
can be equipped with the restriction of the weak$^*$-topology on
the space of all such measures. With respect to this topology, the
map $\iota$ is a homeomorphism of ${\cal M\cal F}$ onto ${\cal
M\cal L}$ which factors to a homeomorphism of the space ${\cal
P\cal M\cal F}$ of projective measured foliations  onto the space
${\cal P\cal M\cal L}$ of \emph{projective measured laminations}, i.e.
the quotient of ${\cal M\cal L}$ under the natural action of
$(0,\infty)$. This homeomorphism maps the space of minimal 
projective measured foliations onto the
space ${\cal M\cal P\cal M\cal L}$ of projective measured geodesic
laminations whose support is a minimal geodesic lamination which
fills up $S$. Since
every minimal geodesic lamination is the support of a transverse
translation invariant measure (compare the expository article
[Bo] for a discussion of this fact and related results),
the image of ${\cal M\cal P\cal M\cal L}$ under the
natural forgetful map $\Pi$ which assigns to a projective measured
geodesic lamination its support equals the set ${\cal B}$. As a
consequence, our above theorem is just a reformulation of the
result of Klarreich provided that the coarse Hausdorff topology on
${\cal B}$ is induced from the weak$^*$-topology on ${\cal M\cal P\cal M
\cal L}$ via the surjective map $\Pi$.

For this it suffices to show that the map $\Pi$ is continuous and
closed. To show continuity, let $(\mu_i) \subset{\cal M\cal P\cal
M \cal L}$ be a sequence of projective measured geodesic laminations. Assume
that $\mu_i\to \mu\in {\cal M\cal P\cal M\cal L}$ in the
weak$^*$-topology, so that the support $\Pi(\mu)$ of $\mu$ is
contained in ${\cal B}$. Since the space of geodesic laminations
equipped with the Hausdorff topology is compact, up to passing to
a subsequence we may assume that the laminations $\Pi(\mu_i)\in
{\cal B}$ converge as $i\to \infty$ in the Hausdorff topology to a
geodesic lamination $\tilde \lambda$. Then $\tilde \lambda$
necessarily contains the support $\Pi(\mu)\in {\cal B}$ of $\mu$
as a sublamination and therefore $\Pi(\mu_i)\to \Pi(\mu)$ in the
coarse Hausdorff topology. Note however that $\tilde \lambda$
may contain isolated leaves which are not contained in the support
of $\mu$ [CEG]. 
Since ${\cal M\cal P\cal M\cal L}$ and
${\cal B}$ are Hausdorff spaces, this shows that the map $\Pi$ is
indeed continuous.

To show that the map $\Pi$ is closed, let $A\subset {\cal M\cal P
\cal M\cal L}$ be a closed set and let $(\mu_i)\subset A$ be a
sequence with the property that $(\Pi(\mu_i))\subset {\cal B}$
converges in the coarse Hausdorff topology to a lamination
$\lambda\in {\cal B}$. Up to passing to a subsequence we may
assume that the geodesic laminations $\Pi(\mu_i)$ converge in the
usual Hausdorff topology to a lamination $\tilde \lambda$
containing $\lambda$ as a sublamination. Since the space of
projective measured laminations is compact, after passing to another
subsequence we may assume that the projective measures
$\mu_i$ converge in the weak$^*$-topology to a projective measure
$\mu$. Then $\mu$ is necessarily supported in $\tilde \lambda$.
Now $\lambda$ fills up $S$ by assumption and therefore every
transverse measure on $\tilde \lambda$ is supported in $\lambda$.
Thus we have $\mu\in {\cal M\cal P\cal M\cal L}$ and, in
particular, $\mu\in A$ since $A\subset {\cal M\cal P\cal M \cal L}$ is
closed. This shows that $\Pi$ is closed and consequently our
theorem is just the main result of [K].

Klarreich's proof of the above theorem relies on Teichm\"uller
theory and the results of Masur and Minsky in [MM1]. In this note
we give a more combinatorial proof which uses train tracks and a
result of Bowditch [B]. We discuss the relation between the
complex of train tracks and the complex of curves in Section 2.
The proof of the theorem is completed in Section 3.

\section{The train track complex}

A \emph{train track} on $S$ is an embedded
1-complex $\tau\subset S$ whose edges
(called \emph{branches}) are smooth arcs with
well-defined tangent vectors at the endpoints. At any vertex
(called a \emph{switch}) the incident edges are mutually tangent.
Through each switch there is a path of class $C^1$
which is embedded
in $\tau$ and contains the switch in its interior. In
particular, the branches which are incident
on a fixed switch are divided into
``incoming'' and ``outgoing'' branches according to their inward
pointing tangent at the switch. Each closed curve component of
$\tau$ has a unique bivalent switch, and all other switches are at
least trivalent.
The complementary regions of the
train track have negative Euler characteristic, which means
that they are different from discs with $0,1$ or
$2$ cusps at the boundary and different from
annuli and once-punctured discs
with no cusps at the boundary.
We always identify train
tracks which are isotopic. 
Train tracks were probably used for the first time 
by Williams [W] to study recurrence properties of
dynamical systems (I am grateful to Greg McShane for
pointing this reference out to me). They became widely known
through the work of Thurston about 
the structure of the mapping class group. A detailed
account on train tracks can be found in [PH] and [M].

A train track is called \emph{generic} if all switches are
at most trivalent.
The train track $\tau$ is called \emph{transversely recurrent} if
every branch $b$ of $\tau$ is intersected by an embedded simple
closed curve $c=c(b)\subset S$ which intersects $\tau$
transversely and is such that $S-\tau-c$ does not contain an
embedded \emph{bigon}, i.e. a disc with two corners at the
boundary.

Recall that a geodesic lamination for a complete hyperbolic
structure of finite volume on $S$ is a compact subset of $S$ which
is foliated into simple geodesics. Particular geodesic laminations
are simple closed geodesics, i.e. laminations which consist of a
single leaf. A geodesic lamination $\lambda$ is called minimal if
each of its half-leaves is dense in $\lambda$. Thus a simple
closed geodesic is a minimal geodesic lamination. A minimal
geodesic lamination with more than one leaf has uncountably many
leaves. Every geodesic lamination $\lambda$ is a
disjoint union of finitely many minimal components and a finite
number of non-compact isolated leaves. An isolated leaf of
$\lambda$ either is an isolated closed geodesic and hence a
minimal component, or it \emph{spirals} about one or two minimal
components ([CEG], [O]).

A geodesic lamination $\lambda$ is \emph{maximal} if all its
complementary components are ideal triangles or once punctured monogons.
A geodesic lamination
is called \emph{complete} if it is maximal and can be approximated
in the Hausdorff topology by simple closed geodesics. Every
minimal geodesic lamination is a sublamination of a complete
geodesic lamination [H]. The space ${\cal CL}$ of complete
geodesic laminations on $S$ equipped with the Hausdorff topology
is compact.

A geodesic lamination or a train track $\lambda$ is \emph{carried}
by a transversely recurrent train track $\tau$ if there is a map
$F:S\to S$ of class $C^1$ which is isotopic to the identity and
maps $\lambda$ to $\tau$ in such a way that the restriction of its
differential $dF$ to every tangent line of $\lambda$ is
non-singular. Note that this makes sense since a train track has
a tangent line everywhere. A train track $\tau$ is called
\emph{complete} if it is generic and transversely recurrent and if
it carries a complete geodesic lamination [H].

A half-branch $\tilde b$ in a generic
train track $\tau$ incident on a switch $v$
is called
\emph{large} if the switch $v$ is trivalent and if
every arc $\rho:(-\epsilon,\epsilon)\to
\tau$ of class $C^1$ which passes through $v$ meets the
interior of $\tilde b$.
A branch $b$ in $\tau$ is called
\emph{large} if each of its two half-branches is
large; in this case $b$ is necessarily incident on two distinct
switches
(for all this, see [PH]).

There is a simple way to modify a complete train track $\tau$ to
another complete train track. Namely, if $e$ is a large branch of
$\tau$ then we can perform a right or left \emph{split} of $\tau$
at $e$ as shown in Figure A below.
The split $\tau^\prime$ of a
train track $\tau$ is carried by $\tau$.
If $\tau$ is complete and if
$\lambda\in {\cal CL}$ is carried by $\tau$, then for every large
branch $e$ of $\tau$ there is a unique choice of a right or left
split of $\tau$ at $e$ with the property that the split track
$\tau^\prime$ carries $\lambda$, and $\tau^\prime$ is complete. In
particular, a complete train track $\tau$ can always be split at
any large branch $e$ to a complete train track $\tau^\prime$;
however there may be a choice of a right or left split at $e$ such
that the resulting track is not complete any more (compare p.120
in [PH]).

\begin{figure}[ht]
\includegraphics{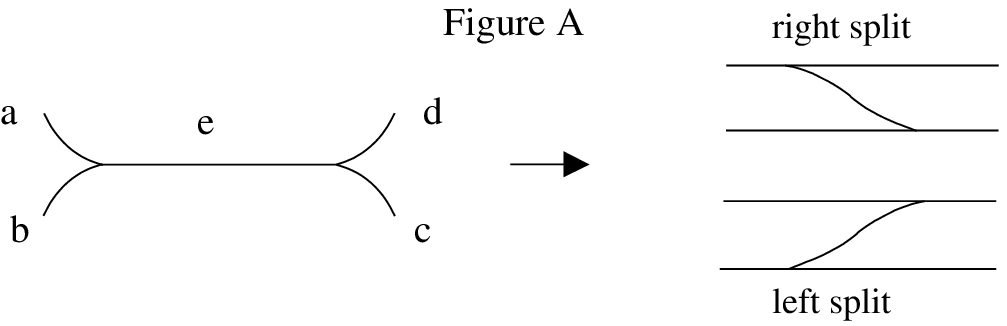}
\end{figure}

Let ${\cal T}T$ be the set of all isotopy classes of complete
train tracks on $S$. We connect two train tracks
$\tau,\tau^\prime$ with a directed edge if $\tau^\prime$ can be
obtained from $\tau$ by a single split at a large branch $e$. This
provides ${\cal T}T$ with the structure of a locally finite
directed metric graph. The \emph{mapping class group} ${\cal
M}_{g,m}$ of all isotopy classes of orientation preserving
homeomorphisms of $S$ acts naturally on ${\cal T}T$ as a group of
simplicial isometries. The following result is shown in [H].

\bigskip

{\bf Theorem 2.1.} {\it The train track complex
${\cal T}T$ is connected, and the action of the mapping
class group on ${\cal T}T$ is proper and cocompact.}

\bigskip

A \emph{transverse measure} on a train track $\tau$ is a
nonnegative weight function $\mu$ on the branches of $\tau$
satisfying the \emph{switch condition}: For every switch $s$ of
$\tau$, the sum of the weights over all incoming branches at $s$
is required to coincide with the sum of the weights over all
outgoing branches at $s$. The set $V(\tau)$ of all transverse
measures on $\tau$ is a closed convex cone in a linear space and
hence topologically it is a closed cell. The train track is called
\emph{recurrent} if it admits a transverse measure which is
positive on every branch. 
A complete train track $\tau$ is recurrent [H].

A transverse measure $\mu$ on $\tau$ is called a \emph{vertex
cycle} [MM1] if $\mu$ spans an extreme ray in $V(\tau)$. Up to
scaling, every vertex cycle $\mu$ is a \emph{counting measure} of
a simple closed curve $c$ which is carried by $\tau$. This means
that for a carrying map $F:c\to \tau$ and every open
branch $b$ of
$\tau$ the $\mu$-weight of $\tau$ equals the number of connected
components of $F^{-1}(b)$. More generally, every integral
transverse measure $\mu$ for $\tau$ defines uniquely a simple
\emph{weighted geodesic multicurve}, i.e. there are simple closed
pairwise disjoint geodesics $c_1,\dots,c_\ell$ and a carrying map
$F:\cup_i c_i\to \tau$ such that $\mu=\sum a_i\nu_i$ where $a_i>0$
is a positive integer and where $\nu_i$ is the counting measure
for $c_i$. We have.

\bigskip

{\bf Lemma 2.2.} {\it Let $c$ be a simple closed
curve which is carried by $\tau$, with carrying
map $F:c\to \tau$. Then $c$ defines a vertex
cycle on $\tau$ if and only if $F(c)$ passes through
every branch of $\tau$ at most twice, with
different orientation.}

{\it Proof:} Let $F:c\to \tau$ be a carrying map for a simple
closed curve $c:S^1\to S$ which defines a vertex cycle $\mu$ for
$\tau$. Assume to the contrary that there is a branch $b$ of
$\tau$ with the property that $F c$ passes through $b$ twice in
the same direction. Then there is a closed nontrivial subarc
$[p,q]\subset S^1$ with nontrivial complement such that $F\circ
c[p,q]$ and $F\circ c[q,p]$ are closed (not necessarily simple)
curves on $\tau$. For a branch $e$ of $\tau$ define $\nu(e)$ to be
the number of components of $(F\circ c[p,q])^{-1}(e)$. Then $\nu$
is a nontrivial nonnegative integral weight function on the
branches of $\tau$ which clearly satisfies the switch condition,
and the same is true for $\mu-\nu$. As a consequence, the
transverse measure $\mu$ can be decomposed into a nontrivial sum
of integral transverse measures which contradicts our assumption
that $\mu$ is a vertex cycle for $\tau$. This shows the first part
of the lemma, and the second part follows in the same way. \qed

\bigskip

In the sequel we mean by a vertex cycle of a complete train track
$\tau$ an \emph{integral} transverse measure on $\tau$ which is
the counting measure of a simple closed curve $c$ on $S$ carried
by $\tau$ and which spans an extreme ray of $V(\tau)$; we also use
the notion vertex cycle for the simple closed curve $c$. As a
consequence of Lemma 2.2 and the fact that the
number of branches of a complete train track on $S$
only depends on the topological type of $S$,
the number of vertex cycles
for a complete train track on $S$ is bounded by a universal
constant (see [MM1]).

Recall that
the \emph{intersection number} $i(\gamma,\delta)$ between two
simple closed geodesics $\gamma,\delta$ equals the minimal number of
intersection points between representatives of
the free homotopy classes of $\gamma,\delta$.
This intersection number extends bilinearly to a pairing for
weighted simple geodesic multicurves on $S$.
The following corollary is immediate from Lemma 2.2.
For its formulation, for a transverse measure $\mu$ on
a train track $\tau$ denote by $\mu(\tau)$ the \emph{total
mass of $\mu$}, i.e. $\mu(\tau)=\sum_b\mu(b)$ where
$b$ runs through the branches of $\tau$. We have.

\bigskip

{\bf Corollary 2.3.} {\it Let $\mu\in V(\tau)$ be an integral
transverse measure on $\tau$ which defines the
weighted simple geodesic
multicurve $c$. Let $\xi$ be any vertex cycle of $\tau$; then
$i(c,\xi)\leq 2\mu(\tau)$.}

{\it Proof:} Let $c$ be any simple closed curve which is carried
by the complete train track $\tau$ and denote by $\mu$ the
counting measure on $\tau$ defined by $c$. Write $n=\mu(\tau)$;
then there is a \emph{trainpath} of length $n$, i.e. a
$C^1$-immersion $\rho:[0,n]\to \tau$ which maps each interval
$[i,i+1]$ onto a branch of $\tau$ and which parametrizes the image
of $c$ under a carrying map $c\to \tau$. We then can deform $\rho$
with a smooth homotopy to a closed curve $\rho^\prime:[0,n]\to S$
which is mapped to $\rho$ by a carrying map and is such that for
each $i\leq n$, $\rho^\prime[i,i+1]$ intersects $\tau$ in at most
one point contained in the interior of the branch $\rho[i,i+1]$.

Now let $\xi$ be any vertex cycle of $\tau$. By Lemma 2.2, $\xi$
can be parametrized as a trainpath $\sigma:[0,s]\to \tau$ which
passes through every branch of $\tau$ at most twice. Then the
number of intersection points between $\sigma$ and $\rho^\prime$
is not bigger than $2n=2\mu(\tau)$. This shows the corollary for
simple closed curves $c$ which are carried by $\tau$. The case of a
general weighted simple geodesic multicurve carried by $\tau$ then
follows from linearity of counting measures and the intersection
form. \qed

\bigskip

Since the distance in ${\cal C}(S)$ between two simple
closed curves $a,c$ is bounded from above by
$2 i(a,c)+1$ [MM1], we obtain from Lemma 2.2 and
Corollary 2.3 the
existence of a number $D>0$ with the property that for
every train track $\tau\in {\cal T}T$, the distance
in ${\cal C}(S)$ between any two vertex cycles
of $\tau$ is at most $D$.

Define a map $\Phi:{\cal T}T\to {\cal C}(S)$ by assigning to a
train track $\tau\in {\cal T}T$ a vertex cycle $\Phi(\tau)$ for
$\tau$. Every such map is roughly ${\cal M}_{g,m}$-equivariant.
Namely, for $\psi\in {\cal M}_{g,m}$ and $\tau\in {\cal T}T$, the
distance between $\Phi(\psi(\tau))$ and $\psi(\Phi(\tau))$ is at
most $D$. Denote by $d$ both the distance on ${\cal T}T$ and on
${\cal C}(S)$. We have.

\bigskip

{\bf Corollary 2.4.} {\it There is a number $C>0$ such that
$d(\Phi(\tau),\Phi(\eta))\leq Cd(\tau,\eta)$ for all
$\tau,\eta\in {\cal T}T$.}

{\it Proof:} Let $\alpha:[0,m]\to {\cal T}T$ be any (simplicial)
geodesic. Then for each $i$, either the train track $\alpha(i+1)$
is obtained from $\alpha(i)$ by a single split or $\alpha(i)$ is
obtained from $\alpha(i+1)$ by a single split. Assume that
$\alpha(i+1)$ is obtained from $\alpha(i)$ by a single split. Then
there is a natural carrying map $F:\alpha(i+1)\to\alpha(i)$. By
Lemma 2.2 and the definition of a split, via this carrying map the
counting measure of a vertex cycle $c$ on $\alpha(i+1)$ defines an
integral transverse measure on $\alpha(i)$ whose total mass is
bounded from above by a universal constant.
Thus by Corollary 2.3, the
intersection number between $c$ and any vertex cycle of
$\alpha(i)$ is bounded from above by a universal constant. Then
the distance in ${\cal C}(S)$ between $c$ and any vertex cycle on
$\alpha(i)$ is uniformly bounded as well [MM1]. This shows the
corollary. \qed

\bigskip

Define a \emph{splitting sequence} in ${\cal T}T$ to be a
(simplicial) map $\alpha:[0,m]\to {\cal T}T$ with the property
that for each $i$ the train track $\alpha(i+1)$ can be obtained
from $\alpha(i)$ by a single split.

We use now a construction of Bowditch [B]. Recall the definition
of the intersection form $i$ on simple geodesic
multicurves. For simple geodesic multicurves 
$\alpha,\beta$ on $S$ with
$i(\alpha,\beta)>0$ and $a>0, r>0$ define
\[L_{a}(\alpha,\beta,r)=\{\gamma \in {\cal C}(S)\mid
\max\{ai(\gamma,\alpha), i(\gamma,\beta)/ai(\alpha,\beta)\}\leq
r\}.\] Our next goal is to 
link the sets $L_a(\alpha,\beta,r)$
to splitting sequences. 
For this recall that a 
\emph{pants decomposition} of $S$ is a collection of $3g-3+m$
pairwise disjoint mutually not freely homotopic simple closed
essential curves on $S$, i.e. these curves are not contractible
and not freely homotopic into a puncture. Let $P=\{\gamma_1,\dots,
\gamma_{3g-3+m}\}$ be a pants decomposition for $S$. Then there is
a special family of complete train tracks with the property that
each pants curve $\gamma_i$ admits a closed neighborhood $A$
diffeomorphic to an annulus and such that $\tau\cap A$ is
diffeomorphic to a \emph{standard twist connector} depicted in
Figure B. Such a train track clearly carries each pants curve from
the pants decomposition $P$; we call it \emph{adapted} to $P$ (see
[PH]). The set of train tracks adapted to a pants decomposition
$P$ is invariant under the action of ${\cal M}_{g,m}$.
\begin{figure}[ht]
\includegraphics{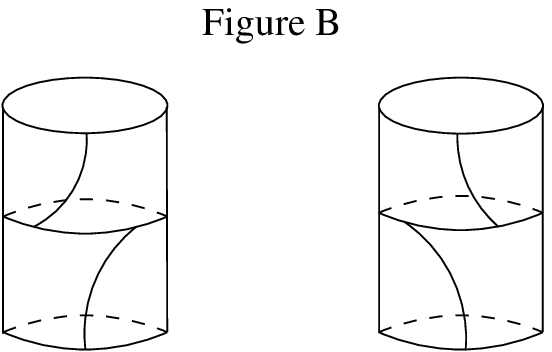}
\end{figure}
We show.

\bigskip

{\bf Lemma 2.5:} {\it 
There is a number $k\geq 1$ with the following
property. Let $\tau_0\in {\cal T}T$ be adapted to a pants
decomposition $P$ of $S$, let
$(\tau_i)_{0\leq i\leq m}\subset {\cal T}T$ be a
splitting sequence issuing from $\tau_0$ and let
$\alpha$ be a simple multicurve consisting of
vertex cycles for $\tau_m$.
Then there is a monotonous
surjective function $\kappa:(0,\infty)\to \{0,\dots,m\}$ such that
$\kappa(s)=0$ for all sufficiently small $s>0, \kappa(s)=m$ for
all sufficiently large $s>0$ and that
for all $s\in (0,\infty)$ there is a vertex cycle
of $\tau_{\kappa(s)}$ which is contained in
$L_s(\alpha,P,k)$.}

{\it Proof:} Let $P$
be a pants decomposition for $S$ and let
$\beta$ be an arbitrary simple multicurve on $S$.
Let $k>1$ and assume that there is a curve $\gamma\in {\cal C}(S)$
with the property that $0<c=i(P,\gamma) i(\gamma,\beta)\leq k
i(P,\beta)$. Write $b=i(P,\gamma)/{i(P,\beta)},
a=i(\gamma,\beta)/c$; then $abi(P,\beta)=1$ and $\max\{a
i(P,\gamma), b i(\beta,\gamma)\}\leq k$. As a consequence, we
have $\gamma\in L(P,\beta,k)$.
Thus for the proof of our lemma we only have to show the
existence of a number $k>0$ with the following property. 
Let $\zeta$ be a train track which is adapted to a
pants decomposition
$P=\{\gamma_1,\dots,\gamma_{3g-3+m}\}$ and let
$\zeta:[0,m]\to {\cal T}T$ be a splitting sequence 
issuing from $\zeta(0)=\zeta$. Let $j>0$
be such that the distance in ${\cal C}(S)$ between every vertex
cycle of $\zeta(j)$ and every vertex cycle of $\zeta(0)$ is at
least $3$. Let $\rho$ be a vertex cycle for $\zeta(m)$;
then there is a vertex cycle
$\alpha(j)$ for $\zeta(j)$
such that \begin{equation}
i(\rho,\alpha(j))\bigl(\sum_{i=1}^{3g-3+m}
i(\alpha(j),\gamma_i)\bigr) \leq k\sum_{i=1}^{3g-3+m}
i(\rho,\gamma_i).\end{equation}

Since $\zeta(0)$ is adapted to the pants decomposition $P$,
every pants curve of $P$ is a vertex cycle
for $\zeta(0)$. Moreover, for each $i\leq 3g-3+m$ there is a
branch $b_i$ of $\zeta(0)$ contained in an annulus $A_i$ about
$\gamma_i$ and such that the counting measure $\nu_i$ for
$\gamma_i$ is the \emph{unique} vertex cycle of $\zeta(0)$ which
gives positive mass to $b_i$. Thus the counting measure $\mu$ of
any simple closed curve $c$ which is carried by $\zeta(0)$ can be
decomposed in a unique way as
$\mu=\mu_0+\sum_{i=1}^{3g-3+m}n_i\nu_i$ where $n_i\geq 0$ and
where $\mu_0$ is an integral transverse measure for $\zeta(0)$
with $\mu_0(b_i)=0$ for all $i$. The intersection number of the
curve $c$ with a pants curve $\gamma_i$ equals the $\mu_0$-weight
of the \emph{large} branch $e_i$ of $\zeta(0)$ contained in the
annulus $A_i$. In particular, the intersection number of $c$ and
$\gamma_i$ coincides with the intersection number of $\gamma_i$
and the simple weighted multicurve $c_0$ defined by the transverse
measure $\mu_0$. Moreover, since the complement of $P$ in $S$ does
not contain any essential closed curve which is not homotopic into
a boundary component or
a cusp, there is a constant $k_0$ only depending on
the topological type of $S$ with the property that
\begin{equation}
\mu_0(\zeta(0))\geq \sum_{i=1}^{3g-3+m} i(c,\gamma_i)
=\sum_{i=1}^{3g-3+m} i(c_0,\gamma_i)\geq \mu_0(\zeta(0))/k_0.
\end{equation}

Consider again the splitting sequence
$\zeta:[0,m]\to {\cal T}T$ and 
let $j\leq m$ be such that the distance in ${\cal C}(S)$ between
every vertex cycle of $\zeta(j)$ and every vertex cycle of
$\zeta(0)$ is at least $3$. Let $\rho$ be a vertex cycle for the
train track $\zeta(m)$. Since $\zeta(m)$ is
carried by $\zeta(j)$, the curve $\rho$ defines a counting measure
$\eta$ on $\zeta(j)$. This counting measure can (perhaps
non-uniquely) be written in the form $\eta=\sum_{i=1}^d a_i \xi_i$
where $\xi_i$ $(i=1,\dots,d)$ are the vertex cycles of $\zeta(j)$
and $a_i\geq 0$ are nonnegative integers. The number $d$ of these
vertex cycles is bounded from above by a universal constant and by
Lemma 2.2, the total mass of each of these vertex cycles $\xi_i$
is bounded from above by a universal constant as well. Therefore,
there is a universal number $q>0$ and there is some $i\leq d$ such
that $a_i\geq \eta(\zeta(j))/q$. After reordering we may assume
that $i=1$. Write $\xi=\xi_1$; Corollary 2.3 shows
that $i(\rho,\xi)\leq
2\eta(\zeta(j))\leq 2qa_1$.

On the other hand, by our assumption on $\zeta(j)$ the distance in
${\cal C}(S)$ between $\xi$ and each of the curves $\gamma_i$ is
at least $3$. Thus $\xi$ is mapped via the carrying map
$\zeta(j)\to \zeta(0)$ to a curve in $\zeta(0)$ which together
with each of the pants curves of $P$ fills up $S$. Then $\xi$
defines a counting measure $\mu$ on $\zeta(0)$, and we have
$\mu=\mu_0+\sum_{i=1}^{3g-3+m}p_i\nu_i$ for some $p_i\geq 0$ with
$\mu_0\not=0$. By inequality (2), the sum of the intersection
numbers between $\xi$ and the curves $\gamma_i$ is contained in
the interval $[\mu_0(\zeta(0))/k_0,\mu_0(\zeta(0))]$. On the other
hand, by the choice of $\xi$ and the fact that the carrying map
$\zeta(j)\to \zeta(0)$ maps the convex cone $V(\zeta(j))$ of
transverse measures on $\zeta(j)$ \emph{linearly} into the convex
cone of transverse measures on $\zeta(0)$, the counting measure
for our curve $\rho$ viewed as a curve which is carried by
$\zeta(0)$ is of the form $a_1\mu +\mu^\prime$; in particular, we
have
\begin{equation}
\sum_i i(\rho,\gamma_i)\geq a_1\mu_0(\zeta(0))/k_0\geq
\eta(\zeta(j))\mu_0(\zeta(0))/qk_0\geq i(\rho,\xi)\mu_0(\zeta(0))/2qk_0.
\end{equation}
As a consequence
of inequalities (2),(3) we have
\begin{equation}
i(\rho,\xi)\sum_i i (\xi,\gamma_i)\leq 2qk_0 \sum_i
i(\rho,\gamma_i).
\end{equation}
This completes the proof of the lemma. \qed
\bigskip

For any metric space $(X,d)$ and any $L\geq 1$, a curve
$\gamma:(a,b)\to X$ is called an \emph{$L$-quasigeodesic} if for
all $a< s<t< b$ we have
\[d(\gamma(s),\gamma(t))/L-L\leq
t-s \leq L d(\gamma(s),\gamma(t))+L.\]

Since ${\cal C}(S)$ is a $\delta$-hyperbolic geodesic metric space
for some $\delta>0$, every $L$-quasigeodesic of finite length is
contained in a uniformly bounded neighborhood of a geodesic in
${\cal C}(S)$. Call a path $\gamma:[0,m]\to {\cal C}(S)$ an
\emph{unparametrized $L$-quasigeodesic} if there is some $s>0$ and
a homeomorphism $\sigma:[0,s]\to [0,m]$ such that the path
$\gamma\circ \sigma:[0,s]\to {\cal C}(S)$ is an $L$-quasigeodesic.
The image of every unparametrized $L$-quasigeodesic in ${\cal
C}(S)$ of finite length is contained in a uniformly bounded
neighborhood of a geodesic.

The following corollary is the key step
toward the investigation of the Gromov boundary of ${\cal C}(S)$.
It was first shown by Masur and Minsky [MM2], with a different
proof.

\bigskip

{\bf Corollary 2.6.} {\it There is a number $Q>0$ such that the image
under $\Phi$ of every splitting sequence in ${\cal T}T$ is an
unparametrized $Q$-quasigeodesic.}

{\it Proof:} Recall the definition of the sets
$L_a(\alpha,\beta,r)$ for $\alpha,\beta\in {\cal C}(S)$.
Bowditch [B] showed that there is a number $r_0>0$ with the
following property. Assume that $\alpha,\beta\in {\cal C}(S)$
\emph{fill up} $S$, i.e. the distance $d(\alpha,\beta)$ between
$\alpha$ and $\beta$ in ${\cal C}(S)$ is at least $3$; then we have.
\begin{enumerate}
\item $L_{a}(\alpha, \beta,r_0)\not=\emptyset$ for all $a>0$.
\item For every $r>0$, $a>0$
the diameter of $L_{a}(\alpha,\beta,r)$ is bounded
from above by a universal constant only depending on $r$.
\item For $r>r_0$ there is a constant $q(r)>0$ with the following
property. For $a>0$ choose some $\gamma(a)\in
L_{a}(\alpha,\beta,r)$; then $\gamma:(0,\infty)\to {\cal C}(S)$ is
an unparametrized $q(r)$-quasigeodesic with
$d(\gamma(s),\alpha)\leq q(r)$ for all sufficiently
large $s>0$
and $d(\gamma(s),\beta)\leq q(r)$ for all sufficiently small
$s>0$.
\end{enumerate}

Let again $\alpha,\beta\in {\cal C}(S)$ be such that
$\alpha,\beta$ fill up $S$.
For $r>r_0$ define
\[L(\alpha,\beta,r)=\cup_{a}
L_{a}(\alpha,\beta,r).\] 
By property 3) above and hyperbolicity of
the complex of curves, there is a number $D(r)>0$ only depending
on $r$ such that $L(\alpha,\beta,r)$ is contained in a tubular
neighborhood of radius $D(r)$ about a geodesic connecting $\alpha$
to $\beta$.

Now let $P$ be any pants decomposition for $S$ 
containing the curve $\alpha$ and assume that
$\gamma\in L_a(P,\beta,r)$ for some $r>0$. Let
$\alpha^\prime$ be a pants curve of $P$ so that
$i(\alpha^\prime,\beta)=\max\{i(\nu,\beta)\mid
\nu\in P\}$; then we have $\gamma\in L_a(\alpha^\prime,\beta,
(3g-3+m)r)$. As a consequence of this,
hyperbolicity of ${\cal C}(S)$ and Lemma 2.5, 
the image under $\Phi$ of the splitting sequence $\zeta$ 
is contained in a uniformly bounded neighborhood 
of any geodesic in ${\cal C}(S)$ connecting
$\alpha$ to $\beta$.
Since this consideration applies to \emph{every}
splitting sequence,
"backtracking" of the
assignment $j\to \Phi(\zeta(j))$ is excluded. From this
the lemma is immediate. \qed

\bigskip

{\bf Remark:} More generally,
the proof of Corollary 2.6 also shows the following. Let
$\zeta,\eta\in  {\cal T}T$ and assume that $\eta$ is
carried by $\zeta$. Let $c$ be any simple closed
curve which is carried by $\eta$; then
$\Phi(\eta)$ is contained in a
uniformly bounded neighborhood of a geodesic arc in ${\cal C}(S)$
connecting $\Phi(\zeta)$ to $c$.

\section{Proof of the theorem}

Fix again a complete hyperbolic metric on $S$ of finite volume.
Recall that a \emph{measured geodesic lamination} on $S$ is a
geodesic lamination equipped with a transverse translation
invariant measure. As in the introduction we equip the space
${\cal M\cal L}$ of measured geodesic laminations with the restriction
of the weak$^*$-topology. The Dirac mass on any simple closed
geodesic $c$ on $S$ defines a measured geodesic lamination. The
intersection of weighted simple geodesic multicurves extends to a
continuous symmetric bilinear form $i$ on ${\cal M\cal L}$ which is
called the \emph{intersection form}. The support of a measured
geodesic lamination $\mu$ for $S$ is minimal and fills up $S$ if
and only if $i(\mu,\nu)>0$ for every measured geodesic lamination
$\nu$ on $S$ whose support does not coincide with the support of
$\mu$. The space ${\cal P\cal M\cal L}$ of \emph{projective measured
laminations} on $S$ is the quotient of ${\cal M\cal L}$ under the
natural action of the multiplicative group $(0,\infty)$; it is
homeomorphic to a sphere of dimension $6g-6+2m-1$ [FLP], in
particular, it is compact. The complex of curves naturally embeds
into ${\cal P\cal M\cal L}$ by assigning to a simple closed geodesic its
projectivized transverse Dirac mass.

Projective measured geodesic laminations can
be used to study infinite sequences in the
complex of curves. Denote again
by $d$ the distance on ${\cal C}(S)$. We have.

\bigskip

{\bf Lemma 3.1.} {\it Let $(c_i)\subset {\cal C}(S)$ be
a sequence which converges in ${\cal P\cal M\cal L}$ to
a projective measured lamination whose support $\lambda_0$
is minimal
and fills up $S$. Let $k>0$ and assume that
$a_i\in {\cal C}(S)$ is such that $d(a_i,c_i)\leq k$; then
up to passing to a subsequence, the sequence $(a_i)$ converges
in ${\cal P\cal M\cal L}$ to a projective measured geodesic lamination
supported in $\lambda_0$.}

{\it Proof:} We use an argument of Luo as explained in the proof
of Proposition 4.6 of [MM1]. Namely, choose a continuous section
$\iota:{\cal P\cal M\cal L}\to {\cal M\cal L}-\{0\}$ of the
projection ${\cal M\cal L}-\{0\}\to {\cal P\cal M\cal L}$. Then
every simple closed geodesic $c$ on $S$ defines a measured
geodesic lamination $\hat c\in \iota({\cal P\cal M\cal L})$. Let
$(c_i)\subset {\cal C}(S)$ be a sequence of simple closed
geodesics. Assume that the sequence $(\hat c_i)$ converges in
$\iota({\cal P\cal M\cal L})$ to a measured geodesic lamination
$\nu_0$ whose support $\lambda_0$ is minimal and fills up $S$.

Let $(a_i)\subset {\cal C}(S)$ be a sequence with $d(a_i,c_i)\leq
k$ for a fixed number $k>0$. By passing to a subsequence we may
assume that $d(c_i,a_i)$ is independent of $i$, i.e. we may assume
that $d(c_i,a_i)=k$ for all $i$. Then for each $i$ there is a
curve $c_i^1\in {\cal C}(S)$ which is disjoint from $c_i$ and such
that $d(c_i^1,a_i) = k-1$. Up to passing to a subsequence, the
sequence $(\hat c_i^1)\subset \iota({\cal P\cal M\cal L})$
converges weakly to a measured geodesic lamination $\nu_1\in
\iota({\cal P\cal M\cal L})$. Since $i(\hat c_i^1,\hat c_i)=0$ for
all $i$, by continuity of the intersection form we have
$i(\nu_0,\nu_1)=0$ and therefore $\nu_1$ is supported in
$\lambda_0$. Proceeding inductively we conclude that up to passing
to a subsequence, the measured laminations $\hat a_i$ defined by
the curves $a_i$ converge in $\iota({\cal P\cal M\cal L})$ to a
measured lamination which is supported in $\lambda_0$. This shows
the lemma. \qed

\bigskip

Consider again the train track complex ${\cal T}T$. For $\tau\in
{\cal T}T$ denote by $A(\tau)\subset {\cal CL}$ the set of all
complete geodesic laminations carried by $\tau$. Then $A(\tau)$ is
open and closed in ${\cal CL}$. Following [H], define a \emph{full
splitting sequence} in ${\cal T}T$ to be a sequence
$\alpha:[0,\infty)\to {\cal T}T$ with the property that for every
$i\geq 0$, the train track $\alpha(i+1)$ is obtained by splitting
$\alpha(i)$ at each of the large branches precisely once. If
$\tau\in {\cal T}T$ is arbitrary and if $\lambda\in {\cal CL}$ is
a complete geodesic lamination which is carried by $\tau$, then
$\lambda$ determines uniquely a full splitting sequence
$\alpha_{\tau,\lambda}$ issuing from $\tau$ by requiring that each
of the train tracks $\alpha_{\tau,\lambda}(i)$ carries $\lambda$,
and $\cap_i A(\alpha_{\tau,\lambda}(i))=\{\lambda\}$ [H]. Recall
the definition of the map $\Phi:{\cal T}T\to {\cal C}(S)$. By
Corollary 2.6, there is a universal number $Q>0$ such that the curve
$i\to \Phi(\alpha_{\tau,\lambda}(i))$ is an unparametrized
$Q$-quasigeodesic in ${\cal C}(S)$. This means that this curve
defines a quasiisometric embedding of the half-line $[0,\infty)$
into ${\cal C}(S)$ if and only if the diameter in ${\cal C}(S)$ of
the set $\Phi(\alpha_{\tau,\lambda}[0,\infty))$ is infinite.

Let ${\cal B}$ be the set of all minimal geodesic laminations
on $S$ which fill up $S$, equipped with the coarse
Hausdorff topology. Recall that ${\cal B}$
is a Hausdorff space. The next statement is immediate from Lemma 3.1.

\bigskip

{\bf Corollary 3.2:} {\it Let $\lambda\in {\cal CL}$
be a complete geodesic lamination which contains
a sublamination $\lambda_0\in {\cal B}$.
Let $\tau\in
{\cal T}T$ be a train track which carries $\lambda$; then the
diameter of the set $\Phi(\alpha_{\tau,\lambda}[0,\infty))\subset
{\cal C}(S)$ is infinite.}

{\it Proof:} Let $\lambda\in {\cal CL}$ be a complete
geodesic lamination which contains a sublamination
$\lambda_0\in {\cal B}$. Assume that $\lambda$ is carried by
a train track
$\tau\in {\cal T}T$.
Denote by
$\alpha_\lambda=\alpha_{\tau,\lambda}$ the full splitting sequence
issuing from $\tau$ which is determined by $\lambda$.
We have to show that the diameter of the set
$\Phi(\alpha_\lambda[0,\infty))$ is infinite. For this recall that
$\cap_i A(\alpha_\lambda(i))=\{\lambda\}$. Since for each $i$ the
curve $\Phi(\alpha_\lambda(i))$ is carried by $\alpha_\lambda(i)$,
the curves $\Phi(\alpha_\lambda(i))$ viewed as projective measured
laminations converge up to passing to a subsequence as
$i\to\infty$ in ${\cal P\cal M\cal L}$ to a projective measured geodesic
lamination which is supported in $\lambda_0$. Thus by Lemma 3.1,
there is no curve $a\in {\cal C}(S)$ with
$d(\Phi(\alpha_\lambda(i)),a)\leq k$ for a fixed number $k>0$ and
all $i$ and hence the diameter in ${\cal C}(S)$ of the set
$\Phi(\alpha_\lambda[0,\infty))$ is indeed infinite.
\qed

\bigskip

As in the introduction, we call a sequence $(c_i)\subset
{\cal C}(S)$ admissible if for a fixed $p\in {\cal C}(S)$ we have
$(c_i,c_j)_p\to \infty$ $(i,j\to \infty)$.
Two admissible sequences $(a_i),(c_i)\subset {\cal C}(S)$ are
equivalent if $(a_i,c_i)_p\to \infty (i\to\infty)$.
The Gromov boundary $\partial {\cal C}(S)$ of ${\cal C}(S)$
is the set of equivalence classes
of admissible sequences in ${\cal C}(S)$. Note
that any quasigeodesic ray in ${\cal C}(S)$ defines
an admissible sequence.
We use Corollary 3.2 to show.

\bigskip

{\bf Lemma 3.3.} {\it There is an
injective map $\Lambda: {\cal
B}\to \partial {\cal C}(S)$.}

{\it Proof:} Fix a pants decomposition $P$ of $S$. Then there is a
finite collection $\tau_1,\dots,\tau_\ell\subset {\cal T}T$ of
train tracks adapted to $P$ with the property that every complete
geodesic lamination $\lambda\in {\cal CL}$ is carried by one of
the tracks $\tau_i$ (see [PH], [H]).
Let ${\cal A}\subset
{\cal CL}$ be the set of all complete geodesic laminations which
contain a sublamination $\lambda_0\in {\cal B}$.
For $\lambda\in {\cal A}$ let $\tau_j$ be a train track
from our collection $\tau_1,\dots,\tau_\ell$ which carries
$\lambda$ and let $\alpha_\lambda:[0,\infty)\to {\cal T}T$ be the
full splitting sequence issuing from $\tau_j$ which is determined
by $\lambda$. By Corollary 2.6 and
Corollary 3.2, there is a universal number $Q>0$ with the
property that
the curve $i\to \Phi(\alpha_\lambda(i))$ is an unparametrized
$Q$-quasigeodesic of infinite diameter. Hence this curve
defines a point $\tilde \Lambda(\lambda)\in \partial {\cal C}(S)$.

There is a natural continuous projection $\pi:{\cal A}\to {\cal
B}$ which maps a lamination $\lambda\in {\cal A}$ to its unique
minimal sublamination $\pi(\lambda)\in {\cal B}$. We claim that
$\tilde \Lambda(\lambda) =\tilde \Lambda(\mu)$ for $\lambda,\mu\in
{\cal A}$ if $\pi(\lambda)=\pi(\mu)=\lambda_0$. For this extend
the map $\Phi$ to the collection of all recurrent train tracks on
$S$ by assigning to such a train track $\sigma$ a vertex cycle
$\Phi(\sigma)$ of $\sigma$. Since the minimal sublamination
$\lambda_0=\pi(\lambda)$ of $\lambda$ fills up $S$ and is carried
by each of the train tracks $\alpha_\lambda(i)$, the image of
$\lambda_0$ under a carrying map $\lambda\to \alpha_\lambda(i)$ is
a \emph{recurrent} subtrack $\hat\alpha_\lambda(i)$ of
$\alpha_\lambda(i)$ which is \emph{large}. This means that
$\hat\alpha_\lambda(i)$ is a train track on $S$ which is a subset
of $\alpha_\lambda(i)$ and whose complementary components do not
contains an essential simple closed curve which is not homotopic
into a puncture. By Lemma 2.2, every vertex cycle for
$\hat\alpha_\lambda(i)$ is also a vertex cycle for
$\alpha_\lambda(i)$ and therefore the distance between
$\Phi(\alpha_\lambda(i))$ and $\Phi(\hat\alpha_\lambda(i))$ is
bounded by a universal constant.

Up to isotopy, the train tracks $\hat\alpha_\lambda(i)$ converge
as $i\to\infty$ in the Hausdorff topology to the lamination
$\lambda_0$ (see [M],[H]). Since $\lambda_0$ is a sublamination of
$\mu$, for every $i>0$ there is a number $j(i)>0$ such that the
train track $\hat\alpha_\lambda(j(i))$ is carried by
$\alpha_\mu(i)$ (see [H]). By the remark following the proof of
Corollary 2.6, this implies that $\Phi(\alpha_\mu(i))$ is contained in
a uniformly bounded neighborhood of
$\Phi(\alpha_\lambda[0,\infty))$. Since $i\geq i_0$ was arbitrary,
the Hausdorff distance between the $Q$-quasigeodesics in ${\cal
C}(S)$ defined by $\lambda,\mu$ is bounded and hence we have
$\tilde \Lambda(\lambda)= \tilde \Lambda(\mu)$ as claimed. Thus
there is a map $\Lambda:{\cal B}\to
\partial {\cal C}(S)$ such that $\tilde \Lambda=\Lambda\circ \pi$.

We claim that the map $\Lambda$ is injective. For this let
$\lambda_0\not=\mu_0\in {\cal B}$ and let $\lambda\in
\pi^{-1}(\lambda_0)\subset {\cal A}, \mu\in \pi^{-1}(\mu_0)\subset
{\cal A}$. By Corollary 2.6 and Corollary 3.2, the image under $\Phi$ of
full splitting sequences $\alpha_\lambda,\alpha_\mu
\in {\cal T}T$ determined by
$\lambda,\mu$ are unparametrized $Q$-quasigeodesics in ${\cal
C}(S)$ of infinite diameter. Thus by the definition of $\Lambda$,
we have
$\Lambda(\lambda_0)=\Lambda(\mu_0)$ if and only if the Hausdorff
distance between $\Phi(\alpha_\lambda [0,\infty))$ and
$\Phi(\alpha_\mu[0,\infty))$ is finite.

Assume to the contrary that this is the case. Then there is a
number $D>0$ and for every $i>0$ there is a number $j(i)>0$ such
that $d(\Phi(\alpha_\lambda(i)),\Phi(\alpha_\mu(j(i))) \leq D$.
Since $d(\Phi(\alpha_\lambda(0)),\Phi(\alpha_\lambda(i)))\to
\infty$ we have $j(i)\to\infty$ $(i\to\infty)$ by Corollary 2.4.
Therefore, up to passing to a subsequence, the curves
$\Phi(\alpha_\lambda(i))$, $\Phi(\alpha_\mu(j(i)))$ viewed as
projective measured geodesic laminations converge as $i\to\infty$
to projective measured geodesic laminations $\nu_0,\nu_1$
supported in $\lambda_0$, $\mu_0$. But $\lambda_0,\mu_0$ fill up
$S$ and do not coincide and hence this contradicts Lemma 3.1.
\qed

\bigskip

The Gromov boundary $\partial {\cal C}(S)$ of ${\cal C}(S)$ admits
a natural Hausdorff topology which can be described as follows.
Extend the Gromov product $(,)_p$ to a product on $\partial {\cal
C}(S)$ by defining $(\xi,\zeta)_p=\sup
\lim\inf_{i,j\to\infty} (x_i,y_j)_p$ where the supremum is taken
over all admissible sequences $(x_i),(y_j)$ 
representing the points $\xi,\zeta$.
We have $(\xi,\zeta)_p=\infty$ if and only if $\xi=\zeta\in
\partial {\cal C}(S)$. A subset $U$ of $\partial {\cal C}(S)$ is
a neighborhood of a point $\xi\in
\partial{\cal C}(S)$ if and only if there is a number
$\epsilon >0$ such that $\{\zeta\in \partial {\cal
C}(S)\mid e^{-(\xi,\zeta)_p}< \epsilon\}\subset U$ (compare [BH]).

We say that a sequence $(c_i)\subset {\cal
C}(S)$ \emph{converges in the coarse Hausdorff topology} to a
lamination $\mu\in {\cal B }$ if every accumulation point
of $(c_i)$ with respect to the
Hausdorff topology contains $\mu$
as a sublamination.
The next lemma
completes the proof of our theorem from the introduction.

\bigskip

{\bf Lemma 3.4:} {\it \begin{enumerate}
\item The map $\Lambda:{\cal B}\to \partial {\cal C}(S)$
is a homeomorphism.
\item
For $\mu\in {\cal B}$,
a sequence $(c_i)\subset {\cal C}(S)$ is admissible and defines
the point $\Lambda(\mu)\in \partial {\cal C}(S)$
if and only if $c_i\to \mu$ in the coarse Hausdorff topology.
\end{enumerate}
}

{\it Proof:} We show first the following. Let $(c_i)\subset
{\cal C}(S)$ be an admissible sequence,
i.e. a sequence with the property that $(c_i,c_j)_p\to
\infty$ $(i,j\to \infty)$.
Then there is some
$\lambda_0\in {\cal B}$ such that $(c_i)$ converges in the coarse
Hausdorff topology to $\lambda_0$.

For this we first claim that there is a
number $b>0$ and an admissible sequence $(a_j)\subset{\cal C}(S)$
which is equivalent to $(c_i)$ (i.e. which satisfies
$(a_i,c_i)_p\to \infty)$ and such that the assignment $j\to a_j$
is a $b$-quasigeodesic in ${\cal C}(S)$.

Namely, let $j>0$ and choose a number $n(j)>j$ such that
$(c_\ell,c_n)_p\geq j$ for all $\ell,n\geq n(j)$. By
hyperbolicity, this means that there is a point $a_j\in {\cal
C}(S)$ with $d(p,a_j)\geq j$ and the property that for $n\geq
n(j)$, every geodesic connecting $c_n$ to $p$ passes through a
neighborhood of the point $a_j$ of uniformly bounded diameter not
depending on $j$. By construction, the sequence $(a_j)\subset{\cal
C}(S)$ is contained in a $b$-quasigeodesic for a number $b>0$ only
depending on the hyperbolicity constant, and this quasigeodesic
defines the same equivalence class as the sequence $(c_j)$. As a
consequence, we may assume without loss of generality that $(c_i)$
is a uniform quasigeodesic. By the considerations in Section 2
we may moreover assume that there is a splitting
sequence $(\tau_j)_{j\geq 0}\subset {\cal T}T$ 
and a strictly increasing
function $\sigma:\mathbb{N}\to \mathbb{N}$ such that
$c_i=\Phi(\tau_{\sigma(i)})$ where $\Phi:{\cal T}T\to
{\cal C}(S)$ assigns to a train track $\tau$ a vertex cycle for
$\tau$.

By Lemma 2.5 there
is a number $k>0$ with the property that
for all $0<i< j$ there is a vertex cycle 
$a_{i,j}\in {\cal C}(S)$ for $\tau_{\sigma(i)}$ 
such that 
\begin{equation}
i(c_0,a_{i,j})i(a_{i,j},c_j)\leq k i(c_0,c_j)\quad \hbox{for}\; 0<i<j.
\end{equation}
Note that this inequality is invariant under multiplication
of the simple closed curve $a_{i,j}$ with an arbitrary positive
weight. Let again $\iota:{\cal P\cal M\cal L}\to {\cal M\cal L}-\{0\}$
be a continuous section and for $j>0$ let $\hat c_j\in
\iota({\cal P\cal M\cal L})$ be a multiple of $c_j$.
By passing to a subsequence we may assume that the
sequence $(\hat c_j)$ converges in the space of measured
geodesic laminations to a measured geodesic lamination
$\mu$.

We claim that the support of $\mu$ is a minimal geodesic
lamination which fills up $S$. For this we argue by contradiction
and we assume otherwise. Then there is a simple closed curve $c$
on $S$ with $i(c,\mu)=0$ (it is possible that the curve $c$ is a
minimal component of the support of $\mu$). Replace the
quasigeodesic $(c_i)$ by an equivalent quasigeodesic, again
denoted by $(c_i)$, which issues from $c=c_0$ and which eventually
coincides with the original quasigeodesic. Such a quasigeodesic
exists by hyperbolicity of ${\cal C}(S)$. 
Since the number of vertex cycles for a fixed
train track is bounded from above by a universal constant,
after passing to a subsequence and using a standard diagonal
argument we may assume that
the curve $a_{i,j}$ is independent of $j>i$; we denote this curve
by $a_i$. Inequality (5) and
continuity of the intersection form then implies that 
$i(c,a_{i})i(a_{i},\hat c_j)\leq k_0 i(c,\mu)= 0$ for
all $i>0$. Since
$d(c,a_{i})\geq d(c,c_i)-k_0$ for all $i$, 
for $i>k_0+2$ the intersection
numbers $i(c,a_{i})$ are bounded from below by a universal
constant and therefore
$i(a_{i},\mu)=0$ for all $i>0$. 
If the support of $\mu$ contains a
simple closed curve component $a$, then this just means that the
set $\{a_i\mid i>0\}\subset {\cal C}(S)$ is contained in the
$k_0+1$-neighborhood of $a$ which is impossible. 
Otherwise $\mu$ has a
minimal component $\mu_0$ which fills a nontrivial bordered
subsurface $S_0$ of $S$, and $i(\mu_0,a)>0$ for every simple
closed curve $a$ in $S$ which is contained in $S_0$ and which is
not freely homotopic into a boundary component or a cusp. Since
$i(a_i,\mu)=0$ by assumption, the curves $a_i$ do not have an
essential intersection with $S_0$ which means that $i(a_i,a)=0$
for every simple closed essential curve $a$ in $S_0$. Again we
deduce that the set $\{a_i\mid i>0\}\subset {\cal C}(S)$ is
bounded. Together we obtain a contradiction which implies that
indeed the support of $\mu$ is a minimal geodesic lamination
$\lambda_0\in {\cal B}$ which fills up $S$.

Let $\lambda_i$ be a complete geodesic lamination
which contains $c_i$ as a minimal component. By passing
to a subsequence we may assume that the laminations $\lambda_i$
converge in the Hausdorff topology to a complete
geodesic lamination $\lambda$. Since the measured laminations
$\hat c_j$ converge in the weak$^*$-topology to $\mu$,
the lamination $\lambda$ necessarily contains $\lambda_0$ as
a sublamination.

Let $\alpha_\lambda$ be a full splitting sequence determined by
$\lambda$. For every $i>0$ the set of complete geodesic
laminations which are carried by $\alpha_\lambda(i)$ is an open
neighborhood of $\lambda$ in ${\cal CL}$. Thus for every $i>0$
there is a number $j(i)>0$ with the property that for every $j\geq
j(i)$ the geodesic $c_{j}$ is carried by $\alpha_\lambda(i)$. From
the remark following Corollary 2.6 we conclude that
$\Phi(\alpha_\lambda(i))$ is contained in a uniformly bounded
neighborhood of any geodesic connecting $c_{j}$ to
$\Phi(\alpha_\lambda(0))$. As a consequence, the image under
$\Phi$ of the full splitting sequence $\alpha_\lambda$ defines the
same point in the Gromov boundary of ${\cal C}(S)$ as $(c_j)$. In
other words, the point in $\partial {\cal C}(S)$ defined by
$(c_j)$ equals $\Lambda(\lambda_0)$ and the map $\Lambda$ is
surjective. Hence by Lemma 3.3, the map $\Lambda$ is a bijection.
Moreover, if $(c_i)\subset {\cal C}(S)$ is any admissible sequence
and if $(c_{i_j})$ is any subsequence with the property that the
curves $c_{i_j}$ converge in the Hausdorff topology to a geodesic
lamination $\lambda$, then $\lambda$ contains a lamination
$\lambda_0\in {\cal B}$ as a minimal component, and $(c_i)$
defines the point $\Lambda(\lambda_0)\in {\cal C}(S)$. In
particular, for every admissible sequence $(c_i)\subset {\cal
C}(S)$ the curves $(c_i)$ converge in the coarse Hausdorff
topology to the lamination $\lambda_0 =\Lambda^{-1}((c_i))\in
{\cal B}$. This shows our above claim.

Let again ${\cal L}$ be the space of all geodesic laminations on
$S$ equipped with the Hausdorff topology. Let ${\cal A}\subset
{\cal L}$ be the set of all laminations containing a minimal
sublamination which fills up $S$. Above we defined a projection
$\pi:{\cal A}\to {\cal B}$. Let $\lambda_0\in {\cal B}$ and let
$L=\pi^{-1}(\lambda_0)\subset {\cal A}$ be the set of all geodesic
laminations which contain $\lambda_0$ as a sublamination. Since
$\lambda_0$ fills up $S$, the set $L$ is finite. We call a subset
$V$ of ${\cal C}(S)\cup {\cal B}$ a \emph{neighborhood} of
$\lambda_0$ \emph{in  the coarse Hausdorff topology of} ${\cal
C}(S)\cup {\cal B}$ if there is a neighborhood $W$ of $L$ in
${\cal L}$ such that $V\supset (W\cap {\cal C}(S))\cup \pi(W\cap
{\cal A})$.

For $\xi\in \partial
{\cal C}(S)$ and $c\in {\cal C}(S)$ write
$(c,\xi)_p=\sup_{(x_i)}\lim\inf_{i\to\infty}(c,x_i)_p$
where the supremum is taken over all admissible sequences
$(x_i)$ defining $\xi$.
A subset $U$ of ${\cal C}(S)\cup \partial{\cal C}(S)$ is called
a neighborhood of $\xi\in \partial {\cal C}(S)$ if there
is some $\epsilon >0$ such that $U$ contains
the set $\{\zeta\in {\cal C}(S)\cup\partial {\cal C}(S)\mid
e^{-(\xi,\zeta)_p}<\epsilon\}$.
In the sequel we identify ${\cal B}$ and $\partial {\cal C}(S)$
with the bijection $\Lambda$. In other words, we view
a point in $\partial {\cal C}(S)$ as a minimal geodesic
lamination which fills up $S$, i.e. we suppress
the map $\Lambda$ in our notation.
To complete the proof of
our lemma it is now enough
to show the following.
A subset $U$ of ${\cal C}(S)\cup
\partial {\cal C}(S)$ is a neighborhood of
$\lambda_0\in {\cal B} =\partial {\cal C}(S)$ if and only if $U$
is a neighborhood of $\lambda_0$ in the coarse Hausdorff topology.

For this let $\lambda_0\in {\cal B}$, let
$L=\pi^{-1}(\lambda_0)\subset{\cal A}$ be the collection of all
geodesic laminations containing $\lambda_0$ as a sublamination and
let $p=\Phi(\tau)$ for a train track $\tau\in {\cal T}T$ which
carries each of the laminations $\lambda\in L$ (see [H] for the
existence of such a train track $\tau$). Let $\epsilon
>0$ and let $U=\{\zeta\in {\cal C}(S)\cup
\partial {\cal C}(S)\mid
e^{-(\lambda_0,\zeta)_p}<\epsilon\}$. Let
$\lambda_1,\dots,\lambda_s\subset L$ be the collection of all
complete geodesic laminations contained in $L$ and for $i\leq s$
let $\alpha_i$ be the full splitting sequence issuing from
$\alpha_i(0)=\tau$ which is determined by $\lambda_i$. By
hyperbolicity and the remark after Corollary 2.6, there is a universal
constant $\chi>0$ with the property that for each $i\leq s,j\geq
0$ every geodesic connecting $p$ to a curve $c\in {\cal C}(S)$ which is
carried by $\alpha_i(j)$ passes through the
$\chi$-neighborhood of $\Phi(\alpha_i(j))$. Since $\Phi(\alpha_i)$
is an unparametrized quasigeodesic which represents the point
$\lambda_0\in \partial {\cal C}(S)$, this implies that there is a
number $j>0$ such that $e^{-(c,\lambda_0)_p}<\epsilon$ and
$e^{-(\mu,\lambda_0)_p}<\epsilon$ for all simple closed curves
$c\in {\cal C}(S)$ and all laminations $\mu\in {\cal B}$ which are
carried by one of the train tracks $\alpha_i(j)$ $(i=1,\dots,s)$.
Since the set of all geodesic laminations which are carried by the
train tracks $\alpha_i(j)$ $(i=1,\dots,s)$ is a neighborhood of
$L$ in ${\cal L}$ with respect to the Hausdorff topology
(see [H]), we
conclude that a neighborhood of $\lambda_0$ in $\partial {\cal
C}(S)\cup {\cal C}(S)$ is a neighborhood of $\lambda_0$ in
the coarse Hausdorff topology as well.

To show that a neighborhood of $\lambda_0\in {\cal B}$ in the
coarse Hausdorff topology contains a set of the form $\{\zeta\in
{\cal C}(S)\cup
\partial {\cal C}(S)\mid e^{-(\lambda_0,\zeta)_p}<\epsilon\}$ we
argue by contradiction. Let again $L=\pi^{-1}(\lambda_0)\subset
{\cal A}$ be the collection of all geodesic laminations containing
$\lambda_0$ as a sublamination. Assume that there is an open
neighborhood $W\subset {\cal L}$
of $L$ in the Hausdorff topology with the property that
$\pi(W\cap {\cal A})\cup (W\cap {\cal C}(S))$
does not
contain a neighborhood of $\lambda_0$ in ${\cal C}(S)\cup
\partial {\cal C}(S)$. Let $(c_i)\subset {\cal C}(S)$ be a
sequence which represents $\lambda_0
\in {\cal B}=\partial {\cal C}(S)$. By our above consideration,
every accumulation point of $(c_i)\subset {\cal L}$
with respect to the Hausdorff
topology is contained in $L$. By our assumption, there is a
sequence $i_j\to\infty$, a sequence $(a_j)\subset {\cal C}(S)$ and
a sequence $R_j\to \infty$ such that $(c_{i_j},a_j)_p\geq R_j$ and
that $a_j\not\in W$. By passing to a subsequence we may assume
that the curves $a_j$ converge in the Hausdorff topology to a
lamination $\zeta\not\in W$. However, since $(c_{i_j},a_j)_p\to
\infty$, the sequence $(a_j)$ is admissible and
equivalent to $(c_j)$ and therefore by our above consideration,
$(a_j)$ converges in the coarse Hausdorff topology to $\lambda_0$.
Then $a_j\in W$ for all sufficiently large $j$ which is a
contradiction. This shows our above claim and completes the proof
of our lemma. \qed

\bigskip

{\bf Acknowledgement:} I am grateful to the referee for pointing
out the reference [MM2] to me and for many other helpful suggestions.

\section{References}

\begin{enumerate}
\item[{[Bo]}] F.~Bonahon, {\em Geodesic laminations
on surfaces}, in
``Laminations and foliations in dynamics, geometry and topology''
(Stony Brook, NY, 1998),  1--37, Contemp. Math., 269,
Amer. Math. Soc., Providence, RI, 2001.
\item[{[B]}] B.~Bowditch, {\em Intersection numbers
and the hyperbolicity of the curve complex}, preprint 2002.
\item[{[BH]}] M.~Bridson, A.~Haefliger, {\sl Metric
spaces of non-positive curvature}, Springer Grund\-leh\-ren 319,
Springer, Berlin 1999.
\item[{[CEG]}] R.~Canary, D.~Epstein, P.~Green,
{\em Notes on notes of Thurston}, in ``Analytical and geometric
aspects of hyperbolic space'', edited by D.~Epstein, London Math.
Soc. Lecture Notes 111, Cambridge University Press, Cambridge 1987.
\item[{[FLP]}] A.~Fathi, F.~Laudenbach, V.~Po\'enaru, {\sl Travaux de
Thurston sur les surfaces,} Ast\'erisque 1991.
\item[{[H]}] U.~Hamenst\"adt, {\em Train tracks
and mapping class groups I}, preprint 2004.
\item[{[Ha]}] W.~J.~Harvey, {\em Boundary structure
of the modular group}, in ``Riemann Surfaces and Related
topics: Proceedings of the 1978 Stony Brook Conference''
edited by I.~Kra and B.~Maskit, Ann. Math. Stud. 97,
Princeton, 1981.
\item[{[I]}] N.~V.~Ivanov, {\em Mapping class groups},
Chapter 12 in ``Handbook of Geometric Topology'', edited by
R.J.~Daverman and R.B.~Sher, Elsevier Science (2002), 523-633.
\item[{[K]}] E.~Klarreich, {\em The boundary at infinity of the curve
complex and the relative Teichm\"uller space,} preprint 1999.
\item[{[L]}] G.~Levitt, {\em Foliations and laminations
on hyperbolic surfaces,} Topology 22 (1983), 119-135.
\item[{[MM1]}] H.~Masur, Y.~Minsky, {\em Geometry of the
complex of curves I: Hyperbolicity}, Invent. Math. 138 (1999),
103-149.
\item[{[MM2]}] H.~Masur, Y.~Minsky, {\em Quasiconvexity in the
curve complex}, preprint 2003.
\item[{[M]}] L.~Mosher, {\em Train track expansions of measured
foliations}, unpublished manuscript.
\item[{[O]}] J.~P.~Otal, {\sl Le Th\'{e}or\`{e}me d'hyperbolisation
pour les vari\'et\'es fibr\'ees de dimension 3}, Ast\'erisque 235,
Soc. Math. Fr. 1996.
\item[{[PH]}] R.~Penner with J.~Harer, {\sl Combinatorics
of train tracks}, Ann. Math. Studies 125, Princeton University
Press, Princeton 1992.
\item[{[W]}] R.F.~Williams, {\em One-dimensional non-wandering sets},
Topology 6 (1967), 473-487.
\end{enumerate}

\end{document}